# Let's Expand Rota's Twelvefold Way For Counting Partitions!

Robert A. Proctor

Department of Mathematics
University of North Carolina
Chapel Hill, North Carolina  27599
rap -at- email.unc.edu

January 5, 2007

*Dedicated to the memory of Gian-Carlo Rota.*

## 1.  Introduction

Let's suppose we have a diamond, an emerald, and a ruby.   How many ways are there to store these gems in small velvet bags?   The bags are identical,  and we can use any number of them.   Unused bags are ignored.   Since the gems can move around within the bags,  we may as well list the gems within each bag in alphabetical order.   Since the bags are identical,  we may as well record the bags in alphabetical order,  according to their first elements.   Using inner  {}'s for the bags and outer  {}'s  to collect the bags,  here are the possibilities:

  { {D}, {E}, {R} },   { {D,E}, {R} },   { {D,R}, {E} },   { {D}, {E,R} },   { {D,E,R} }.

So there are  5  storage configurations.   These are the "partitions" of the set  {D, E, R}.   For  m ≥ 0,  let  B(m)  denote the number of partitions of a set of  m  distinguishable items.   So  B(3) = 5.   The counts  B(0), B(1), B(2), B(3), …  are called the "Bell numbers".   Surprisingly,  these numbers also arise in a Maclaurin series expansion:

$$e^{e^{x}-1} \;=\; \sum_{m \geq 0} B(m)\frac{x^{m}}{m!} \;.$$

There are dozens of variations of this notion of partition.  If autographed baseballs were stored in socks,  the order in which the baseballs were shoved into each sock would matter.   The socks then might be lined up on a shelf,  or thrown into a sack.   One might also consider storing identical un-autographed baseballs.   There would be a limit on how many balls could be stored in one sock.   For each such storage scenario it is natural to ask:   How many configurations can be formed?

If such "grouping" situations are considered one at a time,  a bewildering zoo of enumeration problems arises. Fortunately, in the 1960's Gian-Carlo Rota organized many of the most prominent counting problems of this kind (and their answers) into a nice table with four rows and three columns. Inspired by the Eightfold Way of Buddhism, Joel Spencer suggested naming this table "The Twelvefold Way". The table was publicized by (Rota's student/our advisor) Richard Stanley [**4**, p. 33] and now appears in Doug West's manuscript [**7**, Section 3.4]. Kenneth Bogart added two rows to Rota's table [**1**, p. 115]. We and our Chapel Hill colleague Tom Brylawski (another Rota student) have been independently adding different columns to the Twelvefold Way in our classes. This article publicizes Bogart's two rows and presents the two Chapel Hill columns. Our expanded table has six rows and five columns.

Let's present two more counting problems of this kind. In how many ways can five identical black marbles be assembled into identical plastic sandwich bags scattered around a desktop? Again ignore any unused bags. Here are the possible ways to do this:

$$\{\{\bullet\},\{\bullet\},\{\bullet\},\{\bullet\},\{\bullet\}\}, \quad \{\{\bullet,\bullet\},\{\bullet\},\{\bullet\},\{\bullet\}\}, \quad \{\{\bullet,\bullet\},\{\bullet,\bullet\},\{\bullet\}\},$$

$$\{\{\bullet,\bullet,\bullet\},\{\bullet\},\{\bullet\}\}, \quad \{\{\bullet,\bullet,\bullet\},\{\bullet,\bullet\}\}, \quad \{\{\bullet,\bullet,\bullet,\bullet\},\{\bullet\}\}, \quad \{\{\bullet,\bullet,\bullet,\bullet,\bullet\}\}.$$

More simply, we could just record the number of marbles per bag: $1 + 1 + 1 + 1 + 1$, $2 + 1 + 1 + 1$, $2 + 2 + 1$, $3 + 1 + 1$, $3 + 2$, $4 + 1$, $5$. We have listed all 7 of the "partitions" of the integer 5. For $m \geq 0$, let $p(m)$ denote the number of partitions of the integer $m$. So $p(5) = 7$. One of the most common occurrences of $p(m)$ is as the number of congruence classes in the symmetric group $S_m$. Although the study of $p(m)$ goes back to Euler, there is no known closed-form expression for this quantity. However, Hardy and Ramanujan proved that

$$p(m) \sim (4m\sqrt{3})^{-1} e^{\pi\sqrt{2m/3}} \text{ as } m \to \infty.$$

Many prominent partition counting problems involve a second integer parameter, $b \geq 0$. In the two problems above we could have additionally required that exactly $b$ non-empty bags were to be used. Requiring this with $b = 2$ would have diminished the counts above from 5 to 3 and from 7 to 2. Here's another such problem: a small number $m$ of skiers board a sequence of $b$ large seatless trams. Unused trams stay on the cable and are not ignored. Each skier has a choice of $b$ trams, and so the number of transport configurations of the skiers is $b^m$. This is the number of functions from $\{1,2,\dots,m\}$ to $\{1,2,\dots,b\}$.

The Twelvefold Way puts the counting formulas for many kinds of fundamental combinatorial objects into one framework. These objects include: subsets, multisets, words, permutations, combinations, functions between finite sets, distributions, partitions of sets, and partitions and compositions of integers.  When encountered at an appropriate juncture,  Rota's table can enhance understanding and increase satisfaction for mathematicians and students.  We feel that the most fundamental row in the expanded table is the first of the two new rows.  Giving the two new rows stature equal to the original four rows leads to the adoption of a new conceptual framework:  Rather than partitioning a given collection of items into parts,  we begin with a number of scattered items and then group them into batches.  (This should win us points from "feel good" campus administrators:  we are replacing the old "divisive" approach with a "new age" approach of bringing items together!  :-) )  The traditional terminologies for the original twelve cases were not well suited for unification.  Expansion of the table to thirty cases and the adoption of the grouping viewpoint make accomodation of the traditional terminology even more awkward.  Terminology which facilitates the unified consideration of the thirty cases in this article is presented in Section 2.  (We are not proposing that this terminology be adopted elsewhere.)  As the new viewpoint and terminology are presented in Section 2,  it is explained how each of the  30  empty boxes in Table 1 specifies one set of fundamental combinatorial objects.  Each of these boxes consequently also specifies one enumeration Problem (finding the size of the set) and one Count (that size).  The presentation of the formulaic descriptions of these counts (Answers) in Table 2 is delayed until Section 4 for the sake of *de-emphasis*.

The intended audience for this article begins with combinatorics instructors.  In a limited amount of time,  they must use traditional terminology to teach the most important half dozen or so cases to beginners.  In contrast,  we want to present an "ideal" thirty case version of the Twelvefold Way.  This will be done in a zero-based fashion from the perspective of someone who is familiar with the combinatorics at hand.  We hope that viewing this elementary material from the unifying viewpoints of Tables 1 and 2 will influence the perspectives of combinatorics instructors.  One small direct application of this article could be the adaptation of various small portions of the tables to the course at hand,  to unify and summarize several results at a time.

Years ago we simply wanted to enlarge Rota's table to include a few more interesting cases.  Aesthetic trade-offs soon emerged.  The addition of a row or column can bring along some uninteresting problems or inelegant formulas.  Paul Erdös liked to imagine a heavenly Book in which is written the ideal proof of each theorem in mathematics.  We wondered:  What would The Book's version of Rota's table look like?  It would seem that It must be a fully rectangular table with no footnotes,  which is not too large.  Then we learned of Brylawski's column.   Mixing interactions with Tom and Doug West together with personal aesthetics, delusions of grandeur, and a researcher's curiosity has produced some unanticipated perspectives.   Once the material in Sections 2, 3, 5, and 7 has been presented,  in Section 8 we are in a position to attempt to persuade combinatorialists that it is useful to think of Rota's table as a subtable of Table 2.   We feel that the full framework for the Twelvefold Way is best established by forming the "product" of Brylawski's Column 0 with Bogart's Row A;  doing so produces Table 1.

Our "mother-of-all-twelvefold-way-tables" Table 2 is too large for most courses,  and it would be counterproductive to present too many of the problems from Table 1 at once.   Still, one would hope that most courses could eventually introduce some unifying perspectives rather than just leaving the students with a grab bag of counting fomulas.   The empty-box Table 1 has been presented partly for the convenience of combinatorics instructors.   Despite the fantasy of having formulated the "ideal" table,  stone is the last medium upon which Tables 1 and 2 should be distributed.   Since each instructor will want to abridge or modify these tables differently,  we have posted the Microsoft Word versions of these tables on our website.   Section 2 presents a deliberate *over*-supply of terminology,  from which each instructor may choose a subset to suit his or her pedagogical strategies.   Section 6 contains teaching remarks.

*Monthly* articles are supposed to stimulate and to challenge,  and are meant to be discussed!   To stir the pot a little,  we'll put forward a claim:  *Each of the most fundamental counting quantities somehow appears somewhere in Table 2.*   Note that  2* appears as Answer C4,  for the number  $2^{m-1}$ of "compositions" of  $m \geq 1$  (and not immediately as  $2^n$,  the number of subsets of an  n-set).   And  m!  appears as Answer A2 when  $b = 1$,  for the number of arrangements of  m  books upon one shelf.   But readers who focus on the formulas in Table 2 may miss much of the purpose of this article.   The journey is more important than the destination!   Learning the structure of Table 1 in Section 2 and the relationships amongst the problems that it presents in Section 3 should strengthen the understanding of several fundamental combinatorial concepts.   In Section 4 we advertise some nice-but-underutilized notations when the formulas are presented.

Section 7 contains resource material for student projects involving Neil Sloane's On-Line Encyclopedia of Integer Sequences (OEIS).   In that section we indicate how the extended tabular viewpoint and the writing of this article inspired the filling out and organization of  42  mostly-existing entries in the OEIS.

## 2. The Thirty Assemblage Counting Problems of Table 1

Let $m \geq 1$ and $b \geq 1$. Suppose that we have $m$ marbles scattered around a desktop. We want to count the number of ways in which we can group these marbles into $b$ batches. But first we must specify four attributes of the marbles and of the grouping method: the marbles might all be black or they might have distinct colors; the marbles within a batch might be arranged in a row or just bunched together; the batches might be arranged in a row or scattered about; and the number of marbles in each batch might be unconstrained, be required to be at least one, or be limited to be at most one. Once the four attributes have been specified, each marble configuration which can be so formed is called an assemblage (of a specific kind).

In addition to the grouping viewpoint, there are two other viewpoints for the formation process. Let $S$ be a set of distinct items, or a multiset based upon one type of item. The three viewpoints are: (i) (as above): Start with the items scattered about, and group them into batches, (ii) Split $S$ up into batches, or (iii) Start with some identical empty bins, distribute the items from $S$ into the bins, and then downplay the bins. Although (ii) is the traditional viewpoint for forming partitions, it does not seem to be a natural viewpoint when one wants the items to be ordered within the batches. While (i) will be our primary viewpoint, we will adopt (iii) when it is helpful to have bins available.

As we present the possible attributes of assemblages in detail below in four paragraphs, we explain the structure of Table 1 and introduce associated terminology. Each of the six rows of this table considers certain kinds of assemblages. Each of the five columns of this table specifies a condition on the sizes of the batches within the assemblages. The intersection of a particular row with a particular column determines a set of assemblages of a particular kind. The size of the set of assemblages determined by Row X and Column n is denoted Count Xn. The problem of finding Count Xn is called Problem Xn. The empty boxes of Table 1 specify thirty assemblage counting problems.

When teaching, it is probably best to delay the distribution of Tables 1 and 2 until after the earliest standard counting problems have been studied with traditional terminology. Later, a particular instructor may choose to consider a sequence of problems from some row one week, and then a sequence of problems from some column another week. So in addition to presenting the traditional terminology later in this section, we first supply three kinds of nouns for the 'batch' notion: universal, row-wise, and column-wise. For further flexibility, we also (parenthetically) mention adjective alternatives for the specialized nouns.   Each instructor can use whichever terms seem best for the cases at hand,  and ignore the rest.

Suppose we have  m  items scattered around our desktop.

*How many ways are there to group all of the* m *items into* b *batches?*

We house each batch of items in its own bin.

Are the bins arranged in a row or are they scattered about?   When the bins are arranged in a row,  we say the items are grouped into a <u>sequence</u> of batches (or into <u>serial</u> batches).   When the bins are scattered about,  we say the items are grouped into a <u>collection</u> of batches (or into <u>scattered</u> batches).   Parts I and II of Table 1 present counting problems for assemblages consisting of serial batches and of scattered batches respectively.

Are the items distinguishable or are they identical?   Rows A, B, D, and E of Table 1 present counting problems for assemblages of  m  distinguishable items (into <u>varicolored</u> batches),  while Rows C and F present analogous problems for assemblages of  m  identical items (into <u>monochromatic</u> batches).

Assemblages of  m  items into  b  batches can be produced by distributing the  m  items one at a time into  b  bins.   Suppose the items are distinguishable.   Is the order in which the items arrive in each bin important?   If so,  we say the items are grouped into <u>lists</u> (or <u>filtered</u> batches).   Such assemblages may be visualized by using narrow transparent vertical tubes for the bins.   If the order of arrival is unimportant,  then we say the items are grouped into <u>sets</u> (or <u>unfiltered</u> batches).   Rows A and D of Table 1 present counting problems for assemblages of distinguishable items into lists,  while Rows B and E present analogous problems for assemblages into sets.   Suppose the items are identical.   Then the order of arrival is irrelevant.   A monochromatic batch is a <u>bunch</u>.   Rows C and F present counting problems for assemblages into bunches.

How many items may end up in each batch?   Batches that are required to contain at least one item apiece are called <u>blocks</u> (or <u>non-empty</u> batches).   Batches that are allowed to contain at most one item apiece are called <u>blips</u> (or <u>bitsy</u> batches).   (If the number of items is unconstrained, the batch is <u>freely-sized</u>.)   Columns 1, 2, and 3 of Table 1 present counting problems for assemblages of  m  items into  b  blips,  b  batches,  and  b  blocks respectively.

**Table 1: Thirty Assemblage Counting Problems**

| How many ways to group m items? | Nature of Items | Nature of Batches | (0) $a_i$ Batches Contain i Items | (1) b Blips (# of items $\leq$ 1) | (2) b Batches (any # of items) | (3) b Blocks (# of items $\geq$ 1) | (4) Any Number of Blocks |
|---|---|---|---|---|---|---|---|
| Part I — Sequences of Batches | (A) Disting'ble | Lists | | | | | |
| | (B) Disting'ble | Sets | | | | | |
| | (C) Identical | Bunches | | | | | |
| Part II — Collections of Batches | (D) Disting'ble | Lists | | | | | |
| | (E) Disting'ble | Sets | | | | | |
| | (F) Identical | Bunches | | | | | |

The Twelvefold Way consists of the $4 \times 3$ array of the intersections of Rows B, C, E, and F with Columns 1, 2, and 3. After introducing Rows A and D, Bogart also converted Column 1 into four further rows. He ended up with a $10 \times 2$ table [**1**, p. 115]. So far we have presented an array of $6 \times 3 = 18$ "two-parameter" counting problems.

It's fun to ask students to produce real world models for these counting situations. One popular model for Problems A1-A3 is the grouping of $m$ distinguishable flags using a sequence of $b$ flagpoles. For these problems we prefer another popular model: the grouping of $m$ books using a sequence of $b$ shelves. For Problems D1-D3, the sequence of $b$ shelves can then be replaced by $b$ identical library reshelving carts that are scattered about a large foggy room. The suitability of a model may depend upon the column condition. In Row D, Bogart's model of forming $b$ scattered stacks of books from $m$ books works best for the non-empty lists of Problem D3. People are often used when the items are distinguishable. When working entirely within Column 1, the batch notion of 'blip' should be replaced by the bin notion of 'place'. In Problems A1, B1, D1, and E1, people could be put into chairs arranged in a row or into scattered identical chairs. In Problems C1 and F1, identical light bulbs could be put into lamps arranged in a row or into identical flashlights held in a bag. To emphasize internal disorder, we prefer bags over bins for Rows B and E. To indicate external disorder in Part II, we scatter the identical shelves or bags randomly around the room or the top of a table.

While students may prefer picturesque models, as mathematicians we need succinct descriptions of the possibilities in each of the situations. When the $m$ items are distinguishable, we denote them with the integers $1, 2, \ldots, m$. Then Rows A, B, D, and E consider assemblages of $1, 2, \ldots, m$ into sequences or collections of $b$ lists or sets. Here are examples of one assemblage for each of Rows A, B, D, E when $m = 4$ and $b = 2$:

$$( \, (3,1), (2,4) \, ), \quad ( \, \{2,4\}, \{1,3\} \, ), \quad \{ \, (2,4), (3,1) \, \}, \quad \{ \, \{1,3\}, \{2,4\} \, \} \, .$$

The bins housing the batches of an assemblage appear only subtly in these portrayals, as the inner delimiters () for lists and {} for sets. If the bins are arranged in a row, then () are used as the outer delimiters of each assemblage; if they are scattered about, then {} are used. The assemblages of Problems E3 and B3 have respectively been called "partitions" and "ordered partitions" of the set $\{1,2,\ldots,m\}$. The assemblages of Problems D3 and A3 have at times been called "partitions into lists" and "ordered partitions into lists" respectively.

When the items are identical, each batch may be described by the non-negative integer that counts the items in it. Then the sequences (collections) of bunches of Problem C2 (Problem F2) are described as b-tuples (b-sets) of non-negative integers that sum to m. For collections of bunches, the summands are listed in non-increasing order. Here are example assemblages of these two kinds when m = 9 and b = 4:

$$( 3, 0, 4, 2 ), \quad \{ 4, 2, 2, 1 \} .$$

The summands arising in Problems C3 and F3 must be positive; then they are called "parts" of m. These additive decompositions of the integer m have been called "compositions of m" and "partitions of m" respectively. The sequences of bitsy bunches in Problem C1 become binary sequences of length b that have m 1's in them.

Next we want to consider assemblages which have an unspecified number of batches. If the batches are allowed to be empty, there would be an unlimited number of ways to group m items. But it *does* make sense to count the assemblages of m items into an unspecified number of *blocks*. Here one begins with an unlimited supply of identical bins, the marbles are distributed, and then the empty bins are discarded. Column 4 of Table 1 presents the resulting "one-parameter" counting problems. The jewel and black marble partition counting problems of Section 1 now appear as the m = 3 and m = 5 cases respectively of Problems E4 and F4.

We say a population condition on the batches is <u>homogeneous</u> if it does not refer to the positions of the bins (in Part I), and it is <u>independent</u> if it does not refer to the identities of the items (in Rows A, B, D, E). This article is almost entirely concerned with batch population conditions that are homogeneous and independent. Having already considered simple population conditions for the batches in Columns 1 and 3, following Brylawski we now go all of the way and consider the "finest" homogeneous and independent population conditions: How many assemblages of m items have i items placed in $a_i$ batches for $0 \le i \le m$? Here it must be the case that $\sum_{0 \le i \le m} i a_i = m$, with $a_i \ge 0$. The total number of batches present is $b := \sum_{0 \le i \le m} a_i$. Column 0 of Table 1 presents the resulting multiparameter counting problems. For example, suppose m = 4, $a_2 = 2$, and $a_i = 0$ for $i \ne 2$. Then b = 2 and the possible assemblages for this case of Problem B0 are:

$$( \{1,2\}, \{3,4\} ), \quad ( \{1,3\}, \{2,4\} ), \quad ( \{1,4\}, \{2,3\} ),$$
$$( \{2,3\}, \{1,4\} ), \quad ( \{2,4\}, \{1,3\} ), \quad ( \{3,4\}, \{1,2\} ) .$$

The consideration of completely specified and almost completely unspecified batch population conditions in Columns 0 and 4 has added $6 \times 2 = 12$ counting problems to Table 1.

We struggled with many decisions while endeavoring to develop uniform and comprehensive terminology. Doug West and Rob Donnelly offered valuable feedback . We chose to de-emphasize the distribution viewpoint since it is a dynamic process which gives the bins status equal to that of the items. The bins are nearly invisible in the traditional viewpoint for partitions, and we care only about the final static arrangements produced by distributing. Once the grouping viewpoint was chosen, the noun 'assemblage' was picked as a compromise across the six rows. Within Row D (Row E), a "collection" of lists (sets) is a *set* of lists (sets). Within Row F, a "collection" of bunches is a *multiset* of multisets, each of which is based upon the same one type of item.

Specific terminology may be selected from our system of nomenclature to facilitate the study of cases within one row or one column. For example, when working entirely within Row C, the batches in Columns 1, 2, and 3 can be respectively refered to as bitsy bunches, freely-sized bunches, and non-empty bunches. When working entirely within Column 3, the batches in Rows A, B, and C can be refered to as filtered blocks, unfiltered blocks, or monochromatic blocks. Many synonymous adjective-noun pairs arise, such as 'non-empty list' and 'filtered block'.

We avoid using the adjective 'ordered', since it becomes ambiguous when the noun it modifies is plural. Suppose a visiting algebraist leaves a scrap of paper in your faculty lounge. On it only two words appear: "ordered groups". Is she studying groups whose elements are ordered, or collections of groups which are ordered? "Ordered partitions" is similarly ambiguous to the uninitiated. In fact, it is common to study the lattice formed by the partitions of a fixed set. If we had to choose, we would say that Problem D2 considers unordered carts which hold ordered books. (Rather than unordered collections of ordered carts!)

Doug West suggested 'blip' to denote a batch that can consist of at most one item, since it is close to (binary) 'bit'. Here Rob Donnelly suggested 'places' for the bins. 'Bitsy' comes from the baby talk word 'itsy-bitsy'.

### 3.  Relationships Amongst the Counts of Table 1  (Optional)

Not only will the Editor of The Book take into account the average importance of the problems and the average elegance of the answers within a row or column before accepting the addition of that row or column to Rota's table,  He/She may even consider how nicely the formulas are notated.   But the conceptual framework and structural layout of Table 1 and the relationships amongst its counts are as interesting to us as are the models for these counts and the formulas of Table 2.   This is being emphasized with the presentation of the empty-box Table 1 and by the early placement of this somewhat technical section:   No references to the answers in Table 2 should be made when the relationships in this section are proved.   This section should be skipped by the more casual readers,  since there are only a few references to its content in later sections.

Let  $\mathcal{V}$  be one of the thirty sets of assemblages defined by Table 1.   If the assemblages in  $\mathcal{V}$  were constructed using identical bins placed in a row,  then  let  $\rho$  be the equivalence relation induced within  V  by disregarding the positions of the bins.  (Scatter the bins.)   If the assemblages in  $\mathcal{V}$  were constructed using filtered bins,  then  let  $\sigma$  be the equivalence relation induced within  $\mathcal{V}$  by removing the internal structure from the bins.  (Replace the vertical tubes with bags.)   If the assemblages in  $\mathcal{V}$  were constructed using distinguishable items,  then  let  $\tau$  be the equivalence relation induced within  $\mathcal{V}$  by making the items identical.  (Paint the marbles black.)   If  $\pi \in \{\rho, \sigma, \tau\}$,  saying "apply the symmetry  $\pi$" means "consider the equivalence classes defined by the equivalence relation  $\pi$".   The symmetry  $\sigma$  must be applied before the symmetry  $\tau$  may be applied.

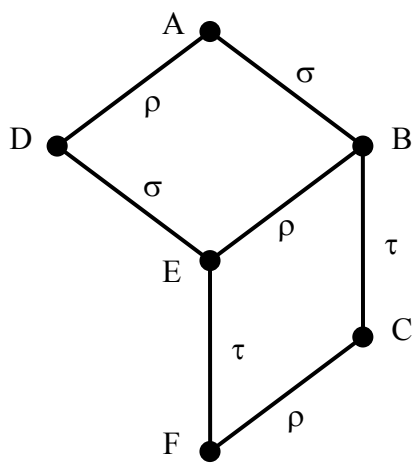

Figure 1

Let's fix a batch population condition and work within its column. The set of assemblages defined by a row below Row A may be viewed as arising from one or more of the sets of assemblages in this column defined by higher rows via the application of one or more of the symmetries $\rho, \sigma, \tau$. Out of the $\binom{6}{2} = 15$ conceivable such relationships amongst the six sets of assemblages, twelve actually arise in this fashion. The poset in Figure 1 depicts the possible applications of the symmetries.

Whenever we are working within Columns 1-3, it is assumed that the two independent parameters $m \geq 1$ and $b \geq 1$ have been fixed. Within Column 4, just $m \geq 1$ has been fixed. When working within Column 0, it is assumed that after $m \geq 1$ is fixed, the $a_i \geq 0$ for $0 \leq i \leq m$ have been chosen such that $\sum_{0 \leq i \leq m} i a_i = m$. Here $b := \sum_{0 \leq i \leq m} a_i$ is a *dependent* variable. We define Columns 1′ and 3′ to respectively consist of the cases in Columns 1 and 3 which may be realized in the context of Column 0. These are the cases which respectively satisfy the natural assumptions $b \geq m$ and $b \leq m$: The counts for all of the other cases of Columns 1 and 3 are 0. (All of the cases of Columns 2 and 4 can realized in the context of Column 0.) Column 0 is the best context for the first result below; three (one) of its parts can be sensibly stated for Column 1′ (3′), but not for Column 1 (3).

An equivalence relation is <u>homogeneous</u> if the equivalence classes it defines have the same size. The notations $\iota$, $\iota!$, and $\binom{m}{\iota}$ are defined at the end of Section 4.

**Quotient Relationships.** *Introduce the equivalence relations $\rho, \sigma, \tau$ within the sets of assemblages from the various rows of Table 1 as indicated in Figure 1. Of the quotient relationships amongst the counts of Table 1 that arise from homogeneous equivalence relations, here are the most important ones:*

*(a) When considering sequences of batches (Part I), in every column the count of assemblages of $m$ identical items into bunches (Row C) can be obtained from the count of assemblages of $m$ distinguishable items into lists (Row A) by dividing by $m!$.*

*(b) When the number $b - a_0$ of non-empty batches is fixed (Columns 0, 1′, and 3′) and the items are distinguishable, the counts of assemblages into collections of $b$ batches (Rows D and E) can be obtained from the counts of the corresponding assemblages into sequences of $b$ batches (Rows A and B) by dividing by $b!/a_0!$.*

*(c)   When the population statistics are completely fixed (Columns 0 and 1´) and the items are distinguishable,  the counts of assemblages into sets with non-empty sizes  ι  (Rows B and E) can be obtained from the counts of the corresponding assemblages into lists (Rows A and D) by dividing by  ι!.*

*(d)   When considering sequences of batches (Part I) in which the population statistics are fixed (Columns 0 and 1´),  the counts of assemblages of  m  identical items into bunches (Row C) whose non-empty sizes are  ι  can be obtained from the counts of the corresponding assemblages of  m  distinguishable items into sets (Row B) by dividing by  $\binom{m}{\iota}$.*

All twelve of the equivalence relations are homogeneous within Column 0.   Parts (a), (b), (c), (d) above respectively describe one, two, two, one of the twelve quotient relationships for counts within Column 0.   Five more quotient relationships arise within Column 0 when whichever as-yet unapplied symmetries are successively applied to group all of the assemblages at hand into one equivalence class,  thereby producing the unity Count F0.   The twelfth quotient relationship within Column 0 is produced by following the symmetry  ρ  from Row A to Row D in Part (b) by the symmetry  σ  from Row D to Row E in Part (c).   All twelve quotient relationships also hold within Column 1´,  but only the Row A to Row C and Row B to Row C relationships of Parts (a) and (d) have both a non-unity divisor and a non-unity quotient.   The three of the twelve quotient relationships which hold within Column 3´ appear in Parts (a) and (b).   Note that Quotient Relationship (a) from Row A to Row C holds for any batch population condition that is independent,  even though the conditions for its "constituent" Quotient Relationships (c) and (d) from Row A to Row B and from Row B to Row C rarely hold.

**Summation Relationships.**   *Of the relationships amongst the counts of Table 1 that arise from forming unions of sets of assemblages of  m  items,  here are the most important ones:*

*(a)   In every row the count of assemblages into an arbitrary number of blocks (Column 4) is obtained by summing the count of assemblages into  b  blocks (Column 3) as  b  runs from  1  to  m.*

*(b)   When considering collections of batches (Part II),  in each row the count of assemblages into  b  batches (Column 2) is obtained by summing the count of assemblages into  d  blocks (Column 3) as  d  runs from  1  to  b.*

*(c)   In every row:  The count of assemblages into  b  batches (Column 2) is obtained from the*

*finest counts of assemblages (Column 0) by summing over all* $(a_0, a_1, \ldots, a_m)$ *such that* $\sum_{0 \le i \le m} i a_i = m$ *with* $a_i \ge 0$ *and* $\sum_{0 \le i \le m} a_i = b$. *The count of assemblages into* $b$ *blocks (Column 3´) is the smaller sum that is obtained by adding the requirement* $a_0 = 0$. *The count of assemblages into an arbitrary number of blocks (Column 4) is then obtained by dropping the requirement* $\sum_{0 \le i \le m} a_i = b$.

Once it is understood how partitions of $m$ can be described using the notation of $(a_0, a_1, \ldots, a_m)$ with $\sum_{0 \le i \le m} i a_i = m$, it is clear that the sums in Summation Relationship (c) for Columns 2-4 respectively run over all partitions of $m$ into: no more than $b$ parts, exactly $b$ parts, and any number of parts.

**Equality Relationships.** *Here are the most important equalities amongst the counts of Table 1 for assemblages of* $m$ *items:*

*(a)  In every row the count of assemblages into* $b$ *blips (Column 1´) is a special case of the finest counts of assemblages (Column 0):  Take* $a_0 = b{-}m$, $a_1 = m$, *and* $a_i = 0$ *for* $2 \le i \le m$.

*(b)  When considering collections of batches (Part II) with* $b \ge m$, *in each row the count of assemblages into* $b$ *batches (Column 2) is equal to the count of assemblages into an unspecified number of blocks (Column 4).*

*(c)  When the batches can contain at most one item (Column 1) and the items are distinguishable,  the counts of assemblages into sets (Rows B and E) are equal to the counts of the corresponding assemblages into lists (Rows A and D).*

## 4.  The Answer Formulas of Table 2

Don't worry … we're not going to end the paper here,  with only the empty Table 1! Now Table 2 will present formulas which answer the problems posed in Table 1.   A few of the formulas are somewhat mysterious and several others are messy sums.   But seventeen of the thirty are products or quotients of products.   Three of the other quantities are famous.

Here's an opportunity to publicize three notations that are not only attractive in their own right,  but whose use also helps to make the most important $3 \times 2$ portion of Table 2 more elegant.   Let $z \in \mathbb{C}$ and $k \ge 0$.   West's [**7**, Section 1.1] is popularizing the notations $z_{(k)} := z(z{-}1)\cdots(z{-}k{+}1)$ and $z^{(k)} := z(z{+}1)\cdots(z{+}k{-}1)$ for falling factorial and rising factorial respectively.  These appeared in C.L. Liu's [**3**, p. 238].  Define $\binom{z}{k} := z_{(k)}/k!$.  We also like Stanley's [**4**, p. 15] notation $\left(\!\binom{z}{k}\!\right) := z^{(k)}/k!$.  Note that $\left(\!\binom{z}{k}\!\right) = \binom{z+k-1}{k}$, since $z^{(k)} = (z+k-1)_{(k)}$.  Let $n \geq 0$,  and suppose that  $k$  items are being selected from an unlimited supply of items of $n$ types.  Here the double parentheses in $\left(\!\binom{n}{k}\!\right)$ indicate that more than one item of a given type may be selected.  (Stanley likes to make note of the intriguing identity $\binom{-n}{k} = (-1)^k \left(\!\binom{n}{k}\!\right)$,  which is the simplest example of the phenomenon of "combinatorial reciprocity".)

Many authors use the notation  $C(b,m)$  to denote the number of combinations of  $m$  items that can be formed from  $b$  distinguishable items.  But it can be seen that this is just another name for Count C1.  More generally,  we already have the name 'Count Xn' available to denote the number of assemblages defined by Problem Xn of Table 1.  We want each entry in Table 2 to be an actual answer,  an evaluatable formula.  Many authors use  $S(m,b)$  to denote the number of partitions of  $[m] := \{1,2,\ldots,m\}$  into  $b$  blocks.  But since we can refer to this number as Count E3,  here the notation  $S(m,b)$  becomes available for a closely related use.

The inclusion-exclusion method can be used to see that Answer E3 is $[\sum_{0 \leq i \leq b}(-1)^i\binom{b}{i}(b-i)^m]/b!$.  Since this formula will not fit into a Table 2 box,  here we give this lengthy expression the name  $S(m,b)$.  The letter "S" is used here because these numbers are known as the "Stirling numbers of the second kind".  We use  $p(m,b)$  to denote the coefficient of  $x^m$  in the formal power series  $x^b \prod_{1 \leq i \leq b}(1 - x^i)^{-1}$.  Direct counting of the terms in the expansion of this formal power series shows that this coefficient is Answer F3.  The letter "p" is used here because Count F3 is the number of partitions of  $m$  into  $b$  parts.  In Section 1 we mentioned the traditional notation  $B(m)$  for Count E4,  the number of partitions of the set  $[m]$.  There we also mentioned the notation  $p(m)$  for Count F4,  the number of partitions of  $m$.  To satisfy our evaluatable goal,  here we redefine  $B(m) := \sum_{1 \leq b \leq m}S(m,b)$  and  $p(m) := \sum_{1 \leq b \leq m}p(m,b)$.  By Summation Relationship (a),  these sums are Answers E4 and F4.

An uncharitable observer could say that the formal power series description for the $p(m,b)$,  $m \geq 0$,  is merely a recasting of the definition for Count F3 into a new language. However,  once Euler had obtained such formal power series descriptions for  $p(m,b)$  and related quantities,  he and his successors (such as Hardy and Ramanujan) used these descriptions to obtain substantive results (such as the asymptotic approximation to  $p(m)$  stated in Section 1). Nowadays symbolic computation packages can be used to quickly find a coefficient in the expansion of such a formal power series,  thereby saving the programming effort needed to count the solutions to $\sum_{0 \leq i \leq m} i a_i = m$ with $a_i \geq 0$ such that $\sum_{0 \leq i \leq m} a_i = b$.

For Column 0, let $\alpha = (a_0, a_1, \ldots, a_m)$ be such that $\sum_{0 \leq i \leq m} i a_i = m$, with $a_i \geq 0$. Set $b := \sum_{0 \leq i \leq m} a_i$. Let $\iota := (1, \ldots, 1, 2, \ldots, 2, \ldots, m, \ldots, m)$ be the $(b-a_0)$-tuple in which there are $a_i$ entries equal to $i$ for each $1 \leq i \leq m$. Define $\alpha! := a_0! a_1! \ldots a_m!$ and $\iota! := 1! \ldots 1! 2! \ldots 2! \ldots m! \ldots m!$. Set $\alpha^+! := a_1! \ldots a_m!$. Define the multinomial coefficients $\binom{b}{\alpha} := b!/\alpha!$ and $\binom{m}{\iota} := m!/\iota!$. Each of the thirty formulas in Table 2 may now be evaluated.

## 5. Favorite Aspects and Miscellaneous Comments

Of the 18 sets of assemblages considered in Columns 1-3, we regard the set $\mathcal{V}$ defined by Problem A2 as being the most fundamental one. These are the assemblages that are modelled by the placement of $m$ books onto a sequence of $b$ shelves. Each of the other sets of assemblages in Column 2 arises as a set $\mathcal{W}$ of equivalence classes formed by the application of one or more of the symmetries $\rho, \sigma, \tau$ within $\mathcal{V}$. Each of the counts in Columns 1 and 3 in a given row then measures a subset of the set $\mathcal{W}$ in that row. Therefore $b^{(m)}$ is at least as large as any of the other entries in Columns 1-3.

The most difficult quantity to get our hands on, $p(m)$, appears as the last entry, Answer F4, of Table 2. All three symmetries have been applied and the number of blocks is open-ended. The simplest quantity, 1, appears as the bottom left corner entry, Answer F0. Given this answer for Problem F0 and the comment near the end of Section 3 (and ignoring Answers F2-F4), it is particularly easy to confirm the three parts of Summation Relationship (c) within Row F. Two of the least interesting entries in Table 2 are Answers B4 and D4. Nonetheless, the reference material in the On-Line Encyclopedia of Integer Sequences for the sequences formed by these counts is respectively approximately 150 and 70 lines long [**6**, A000670, A000262]. Count B4 is the number of rankings of $m$ teams in a league if ties are possible. Count D4 is the number of stackings of $m$ books into scattered piles. Within Column 0 suppose that $a_k = b$ for some $1 \leq k \leq m$, with $a_i = 0$ for $i \neq k$. Then $m = bk$, and $(bk)!/[b!(k!)^b]$ is the answer to the following manifestation of Problem E0: How many ways are there to form $b$ teams of $k$ people apiece from a pool of $bk$ people?

Stanley emphasized [**4**, p. 31] the function enumeration interpretation of the Twelvefold Way: The functions $f$ from $[m]$ to $[b]$ can be modelled by distributing marbles numbered $1, 2, \ldots, m$ into bins in positions $1, 2, \ldots, b$, as in Problem B2. Then Counts B1 and B3 enumerate such functions that are respectively injective and surjective. What if one wishes to regard all of the elements of the domain [m] as being equivalent to each other? And/or all of the elements of the codomain [b]? These equivalences induce the symmetries $\tau$ and $\rho$ within each of the original sets of functions considered in Row B. The notions of injective and surjective still make sense for the resulting equivalence classes of functions. The counts of these equivalence classes appear at the intersections of Rows C, E, and F with Columns 1, 2, and 3. Any function f: [m] → [b] is determined by its coimage (or marble distribution) $(f^{-1}(1), f^{-1}(2), \ldots, f^{-1}(b))$. This has the ( {}, {}, …, {} ) form of an ordered partition in which empty batches are allowed. Bogart noted that the ( (), (), …, () ) assemblages of Problem A2 into sequences of filtered batches analogously determine "ordered distributions". In this article this usage of 'ordered' is replaced by 'filtered':

**Definitions.** *A <u>filtered distribution</u> (<u>function</u>) from a set* M *to a set* B *consists of an assignment of a list of elements (subset) of* M *to each element of* B *such that each element of* M *appears in exactly one of the lists (subsets).*

Counts A1 and A3 respectively enumerate injective and surjective filtered distributions. There are $b^m$ functions from [m] to [b] and $b^{(m)}$ filtered distributions from [m] to [b].

Columns 1′ and 3′ arose from the imposition of the assumptions of $b \geq m$ and $b \leq m$ on Columns 1 and 3 respectively. Here the simpler answer 1 arises for Problems D1′, E1′, F1′, and the quantity $\binom{m-1}{b-1}$ appearing in Answers A3′, C3′, and D3′ may now be written more suggestively as $\left(\binom{b}{m-b}\right)$.

It's instructive to trace the diminishment of Count A0 to Count F0 using Figure 1. Both impositions of $\sigma$ lead to a decrease by a factor of $\iota!$. The two higher impositions of $\rho$ cause a decrease by $b!/a_0!$, but the lowest imposition of $\rho$ causes a decrease by only $[b!/a_0!][1/\alpha^+!]$. The higher imposition of $\tau$ causes a decrease by $\binom{m}{\iota}$, but the lower imposition of $\tau$ causes a decrease by only $[1/\alpha^+!]\binom{m}{\iota}$. Therefore Count A0 is written $\{[b!/a_0!][1/\alpha^+!]\}\binom{m}{\iota}[\iota!]$ via Row C, but is written $[b!/a_0!]\{[1/\alpha^+!]\binom{m}{\iota}\}[\iota!]$ via Row E. The complete generality of Quotient Relationship (a) comes from the disappearance of $\iota$ in the product identity $\binom{m}{\iota}[\iota!] = m!$ for every case of Column 0. There this product $\binom{m}{\iota}[\iota!]$ arose from the composition of the much more special Quotient Relationships (c) and (d).

In Part I one could leave the setting of this paper and consider non-homogenous batch population conditions.  Counts even finer than those in Column 0 arise when the number of items in each batch is specified.  Let $m_1, m_2, \ldots, m_b \geq 0$ be such that $m_1 + m_2 + \ldots + m_b = m$.  For $1 \leq i \leq b$, require that there be $m_i$ items in the $i$th batch.  A Column 00 could be defined in Part I to consider such assemblages.  The three analogous counts in Column 0 would be obtainable from these three finer Column 00 counts by multiplying each by $\binom{b}{\alpha}$, since the batches of those sizes could now be placed in that many new orders.  Define $\mu := (m_1, m_2, \ldots, m_b)$.  Then $\binom{m}{\mu}$ would be a more familiar way to denote the answer $\binom{m}{t}$ to Problem B00:  via the multinomial theorem, these assemblages form a familiar subset of the $b^m$ possibilities for Count B2.

Is it possible to make sense of our counting problems if $m = 0$ and/or $b = 0$?  Yes! However, several problems arise with the formulas presented in Table 2.  First consider $m = 0$ for Columns 0-3.  When $b = 0$, there is 1 way to group zero items:  the empty sequence or the empty collection of batches.  When $b \geq 1$, the counts are again 1 for Columns 1 and 2.  They are 1 for Column 0 if and only if $a_0 = b$; otherwise they are 0.  They are 0 for Column 3. Next consider $m \geq 1$ and $b = 0$.  Here the counts in Columns 0-3 are 0.  The counts in Column 4 are 1 when $m = 0$.  Some of these facts are used in Section 7.

## 6.  Teaching Remarks

Some readers may have noticed that the familiar counting quantities $b_{(m)}$, $b^m$, $\binom{b}{m}$, and $\binom{b+m-1}{m}$ are not portrayed in Table 2 as counting selections:  Isn't $\binom{b}{m}$ supposed to be the number of ways of selecting $m$ out of $b$ distinguishable items?  Instead, the familiar formula $\binom{b}{m}$ appears in Table 2 as the number of ways of assembling $m$ identical items into a sequence of $b$ blips!

It might be best to greatly delay the distribution of Table 1 and to distribute Table 2 only at the end of the semester.  Even then, some instructors may want to delete some portions of the overly large Tables 1 and 2.

We begin our course by deriving the usual four formulas for counting the ways to select $k$ items from a supply of items of $n$ types:  the order of selection may or may not matter, and selecting more than one item of a type may or may not be allowed.  The $\binom{n}{k}$ count is obtained using an equivalence relation, and the $\binom{n+k-1}{k}$ count is then obtained using "dots and bars".

The multinomial coefficient formula for counting anagrams is derived next.  The students are then ready to be assigned a sequence of the five Row C problems for distributing $m$ identical bolts into identical bins placed in a row.   They can use two of the four selection formulas to answer Problems C1 and C2 once they understand that any problem in Row C models selections (of bin positions) with the order of selection being unimportant  (m  black marbles).   (Rota liked to note that Problem C1 (C2) arises in physics as counting the selection of distinct energy levels by indistinguishable fermions (bosons).)   They can answer Problems C3 and C4 using dots and bars.   Students can use the other two of the four selection formulas to answer Problems B1 and B2 once they understand that any problem in Row B models selections (of bin positions) where the order of selection  (m  numbered marbles) is important.   Answer A2 can be obtained from Answer C2.  The Problems C2-B2-A2 segment of Table 1 appears nicely as the three part Exercise 3.45 of Brualdi's text [**2**],  where one is asked to place  20  different books on  5  shelves:   Part (a) cares only about the number of books per shelf,  and Part (b) cares only about which books end up on which shelf.   Another triptych of easy bookshelf problems can be formed from the A2-A3-A4 segment of Table 1.

After inclusion-exclusion and generating functions are studied,  Answers B3 and F3 can be obtained in lecture.   Now the students can be assigned a sequence of six Column 3 problems for arrangements of  m  books on  b  non-empty shelves or carts in the spirit of [**2**, Exercise 3.45].   Following this,  Table 1 could be distributed as a display of the "product" of the six problems of Column 3 with the five problems of Row C.

Textbooks give little or no attention to the most fundamental enumeration technique, counting-by-listing.   And too often we fail to convey the nature of mathematical research to our students.   Recently the core  $6 \times 3$  portion of Table 1 was distributed early in the semester for the following project:   Student teams were asked to empirically find the counts for four specified  (m,b)-parameter values for each of the eighteen core marble grouping problems.   Each of six two-student teams was assigned three of the problems.   The relationships of the results found in lecture early in the semester to these problems were not known to the students.   Only brief definitions of the terminology appearing in Table 1 appeared on this handout (and only on this handout).   So this was truly a research project!   Each team had to develop a standard form for each of its three assigned kinds of assemblages,  plus a systematic way of listing the possible standard forms.   These  $6 \times 3 \times 4 = 18 \times 4 = 72$  empirical counts were collected and distributed.   (Whether correct or not!)   Later in the semester a scrambled list of the formulaic answers for these  18  problems was also distibuted,  and teams were asked to match the answers with the problems using the published data as a partial guide. The Word documents for these two projects are available on this author's website.

The formula $b^{(m)}$ for Answer A2 can be obtained directly. Then the method of Quotient Relationship (a) can be used to re-obtain $\left(\binom{b}{m}\right) = b^{(m)}/m!$ for Answer C2. This derivation can be related to the quotient argument for selections which produced $\binom{b}{m} = b_{(m)}/m!$: Although that method would first be recast to be a manifestation of Quotient Relationship (d), it could be further recast to be a manifestation of Quotient Relationship (a). This would permit the simultaneous derivations of Answers C1 and C2 from Answers A1 and A2 respectively. In any case, students should be asked why Answer C2 is not Answer B2 divided by $m!$, while Answer C1 *is* Answer B1 divided by m!!

If Column 0 is too much for your students to swallow, you might introduce Column 7 instead. It is defined in the fourth paragraph of the next section, and it was exemplified within Row B with $m = 4$ near the end of Section 2.

## 7. Fibonacci and Catalan Hunting, the OEIS, and Student Projects

Stanley gave 9 combinatorial interpretations of the Fibonacci numbers in Exercise 1.14 of [**4**] and 66 combinatorial interpretations of the Catalan numbers in Exercise 6.19 of [**5**]. If a combinatorialist accepts our claim that the most fundamental combinatorial quantities appear somewhere in Table 2, then two questions he or she may soon ask are: How about the Fibonacci numbers? Or the Catalan numbers? Do either of these famous sequences of combinatorial counts arise in Table 2? In this section we also advertise the On-Line Encyclopedia of Integer Sequences, describe the mutually beneficial interactions between the OEIS and the tabular viewpoint, and point out some source material for student projects.

We computed the counts for each of the 18 cases displayed in Columns 1-3 for $0 \le m \le 8$ and $0 \le b \le 8$ and found that the Fibonacci and the Catalan numbers did not begin to appear as subsequences within the resulting eighteen $9 \times 9$ tables. Many Column 0 counts were also examined. Since the Column 4 cases are indexed by only the parameter m, each of these six cases immediately produces a sequence of counts. We soon realized that the batch population criteria of Column 4 could be modified or analogized in various ways to produce other counting sequences that are each indexed by one parameter, either $m \ge 0$ or $b \ge 0$. Six new batch population criteria will be presented in the additional Columns 5-10 below, giving us a total of  7  "one-parameter" columns.   To continue our search for the two famous sequences,
we computed the first nine terms of each of the  $6 \times 7 = 42$  sequences arising from Columns 4-
10 using analogues of Summation Relationships (c).   We then searched for the Fibonacci
numbers  13 and 21  and the Catalan numbers  132 and 429.

Recall that the population condition that enabled the creation of Column 4 was the
requirement that each batch contain at least one item,  i.e.,  that each batch is a block.   Now
Column 5 is defined to consider assemblages of  $m \geq 0$  items in which each batch/block contains
at least  2  items.   Then the assemblages counted in Problem C5 are the compositions of  m  in
which each part is at least  2.   By [**4**, Exercise 1.14(b)],  the number of these when  $m \geq 1$  is the
Fibonacci number  $F_{m-1}$  (where  $F_0 = 0$  and  $F_1 = 1$).   A further Column 6 considers
assemblages of  $m \geq 0$  items into blocks that each contain *at most*  2  items.   Then the
assemblages counted in Problem C6 are the compositions of  m  in which each part is  1  or  2.
By [**4**, Exercise 1.14(c)],  the number of these is the Fibonacci number  $F_{m+1}$.   To summarize,
Columns 4, 5, and 6 consider assemblages of  $m \geq 0$  items into a varying number (possibly  0)
of batches of sizes from the size sets  {1,2,3,…},  {2,3,4,…},  and  {1,2}  respectively.

After seeing an early version of Table 2 with only Columns 1-4,  Doug West suggested
that we also consider assemblages in which exactly two items are placed in each of the bins.
Starting with  $b \geq 0$,  here one takes   m := 2b.   The sizes of the batches are coming from the size
set  {2};  we introduce Column 7 to consider this population condition.   The assemblages
counted in Problem E7 may be viewed as configurations of  b  non-incident chords that connect
2b  points on a circle which are named  1, 2, …, 2b.   If one then further restricts one's attention
to the special such configurations in which no two chords intersect,  the special subcount that
arises [**5**, Exercise 6.19(n)] is the Catalan number  $\binom{2b}{b}/(b+1)$.   However,  this recipe is not an
independent population condition,  since it refers to the labels on the marbles.   When  b = 2m,
the famous lattice path realization [**5**, Exercise 6.19(h)] of the Catalan numbers  $\binom{2m}{m}/(m+1)$  can
be obtained within Problem C1 by imposing an additional non-homogeneous population
condition:  For  $1 \leq p \leq b$,  do not allow the number of marbles placed in the first  p  blips to
exceed  [p/2].

Seeking an appearance of the Catalan numbers via a batch population condition that is
both homogeneous and independent,  we realized that it is also possible to define counting
sequences which are indexed only by the number-of-batches parameter  b.   This can be done by putting a bound upon the number of items that each of the batches may contain and then allowing the total number of items to vary. In particular, the next column is dual to Column 4 in a certain sense: the population condition that enables its creation is the requirement that each batch contain at most one item, i.e., that each batch is a blip. Given $b \geq 0$, Column 8 considers assemblages of from $0$ to $b$ items into $b$ batches which each contain at most $1$ item. When the items are distinguishable, they are to be numbered $1, 2, 3, \ldots$. The resulting Counts A8, B8, …, F8 for $b \geq 0$ are respectively $\sum_{0 \leq m \leq b} b_{(m)}$, $\sum_{0 \leq m \leq b} b_{(m)}$, $2^b$, $b+1$, $b+1$, $b+1$. Given $b \geq 0$, Column 9 considers assemblages of from $0$ to $2b$ items into $b$ batches which each contain at most $2$ items. Counts F9 for $b \geq 0$ are the triangular numbers $\binom{b+2}{2}$. Given $b \geq 0$, Column 10 considers assemblages of from $0$ to $2b$ items into $b$ batches which each contain $1$ or $2$ items. To summarize, Columns 8, 9, 10 consider assemblages into $b \geq 0$ batches of a varying number (possibly $0$) of items, with the batch sizes respectively coming from the size sets $\{0,1\}$, $\{0,1,2\}$, and $\{1,2\}$.

The Catalan numbers 132 and 429 did not appear in Columns 4-6 for $0 \leq m \leq 8$ or in Columns 7-10 for $0 \leq b \leq 8$. The problem of producing the sequence of Catalan numbers with a batch population condition that is both homogeneous and independent remains open.

As a consolation prize, we did end up with the first nine terms for each of the $6 \times 7 = 42$ enumeration sequences on our hard drive. How many (if any) of these sequences were "new"? The OEIS has become an invaluable tool for researchers. It was easy to paste each of these $42$ sequence starts into the OEIS search box to check out "whazzup?". It would have been daunting to find references for these sequences using human or library resources. To some extent, the viewpoint of this paper has "returned the favor" to the OEIS: Three of these sequences turned out to be new to the OEIS, and were subsequently entered into the OEIS: Counts E9, F9, and B10. For the other $39$ counting sequences, their OEIS entries could be used to give some indication of how well known they had been, and whether their earlier occurences were as partition counts. The OEIS reference A-numbers for the $42$ sequences defined by Columns 4-10 are displayed in Table 3. As it turned out, many of the existing entries for the "known" $39$ sequences contained no indication that the sequence at hand was related to partition counting problems. Table 3 and its supporting definitions now appear on the last screen of the entry for "partitions" in the OEIS index [**6**, …/Sindx_Par.html]; that material organizes and interrelates these $42$ sequences as partition counting sequences.

|       | Col 4  | Col 5  | Col 6  | Col 7  | Col 8  | Col 9  | Col 10 |
|-------|--------|--------|--------|--------|--------|--------|--------|
| Row A | 002866 | 052554 | 005442 | 010050 | 000522 | 082765 | 099022 |
| Row B | 000670 | 032032 | 080599 | 000680 | 000522 | 003011 | 105749 |
| Row C | 011782 | 000045 | 000045 | 000012 | 000079 | 000244 | 000079 |
| Row D | 000262 | 052845 | 047974 | 001813 | 000027 | 105747 | 001517 |
| Row E | 000110 | 000296 | 000085 | 001147 | 000027 | 105748 | 001515 |
| Row F | 000041 | 002865 | 008619 | 000012 | 000027 | 000217 | 000027 |

**Table 3**

The OEIS entries for these 42 sequences can serve as sources for special problems or student projects. Providing proofs and exposition of the encyclopedia's terse material for just about any one of these sequences could constitute an end-of-semester project. Many of the entries give summation expressions, recurrences, and/or generating functions for the sequence at hand. Since some of the entries did not include combinatorial descriptions that could immediately be seen to be equivalent to our assemblage descriptions, to know that we had true matches it was necessary for us to supply proofs that our counts satisfied the given OEIS expressions. For example, the OEIS sequence A047974 matched our Count D6 sequence through $m = 8$. That OEIS entry satisfied the recurrence $f_m = f_{m-1} + 2(m-1)f_{m-2}$. Our Count D6 is the number of assemblages of $1, 2, \ldots, m$ into a collection of lists of length 1 or 2. Given such a assemblage, erase the entry $m$. Doing this when $m$ is in its own list produces the $f_{m-1}$ such assemblages of $[m-1]$. If $m$ is not in its own list, also erase the other number $k$ in the same pair and decrement $k+1, k+2, \ldots, m-1$ by 1 each. Since $1 \leq k \leq m-1$ and since $k$ can be to the right or to the left of $m$ in that pair, the number of such assemblages of $[m-2]$ that is produced is $2(m-1)f_{m-2}$. So our D6 sequence can be seen to satisfy the same recurrence. For a more challenging exercise, consider the OEIS sequence A003011, which matched our Count B9 sequence through $b = 8$. That OEIS entry satisfies the recurrence

$$bf_b = (2b^3 - b^2 + b + 1)f_{b-1} + (-3b^3 + 4b^2 + 2b - 3)f_{b-2} + (b^3 - 2b^2 - b + 2)f_{b-3} !$$

Fortunately, this entry of the OEIS also gave the summation expression $\sum_{0 \leq k \leq 2b} k!U_k$, where $U_k$ is the coefficient of $x^k$ in the expansion of $(1 + x + \frac{1}{2}x^2)^b$. Our Count B9 enumerates the assemblages of $1, 2, \ldots, m$ into a sequence of $b$ sets, each with no more than 2 elements, as $m$ ranges from 0 to $2b$. Since it is not hard to see that $m!U_m$ is the number of such assemblages of $1, 2, \ldots, m$ for $0 \leq m \leq 2b$, we can obtain the same sum. The Column 10 analogue of Summation Relationships (c) of Section 3 was used to prove complete matches with the OEIS sequences which matched three Column 10 count sequences through $b = 8$: Summation expressions for our Counts A10, D10, and E10 can be obtained by summing the corresponding Column 0 count over the possible $\alpha$'s: $a_2 = i$ and $a_1 = b - i$ for $0 \le i \le b$. The resulting sums agree with those given in the OEIS for those three sequences.

To take into account further possible marble multiplicities, Table 1 could be extended in another fashion. For example, suppose we have $m = 2n$ marbles. There could be $n$ black marbles and $n$ white marbles, or there could be $n$ distinct colors for the marbles which occur twice each. At this point it might be best to move to a three dimensional table: The columns would remain the same, but only Rows A, B, D, and E would be retained. The top layer would consider the distinct color case and the bottom layer would consider the case where all of the marbles are black. Intermediate layers would consider additional multiplicity possibilities such as the two just mentioned. The question of whether any of the resulting count sequences already appear in the OEIS may be a suitable student research project.

There are ways by which two-parameter count tables can be converted into sequences for entry into the OEIS. We have not checked to see which of the 18 two-parameter count tables defined in Columns 1-3 already appear in the OEIS in some such fashion.

## 8. Concluding Remarks

Building upon Section 3, here are some remarks which support the augmentation of the Twelvefold Way with Rows A and D and Columns 0 and 4.

We like to think of Table 1 as a "product" of Column 0 with Row A. Each of the 30 counts can be expressed in terms of the Counts A0 as follows: For a given population condition, Counts B0, C0, D0, E0, and F0 can be found from Count A0 using the quotient relationships. Then in Row X, the Count Xn can be found by applying Equality Relationship (a) and Summation Relationship (c) to the relevant cases of the six Column 0 counts, using the Row A cases as a guide. Row A has the nicest set of answers amongst the rows (along with Row C). Counting when symmetry is present is generally more difficult. Since the assemblages of Row A have had none of the symmetries $\rho, \sigma, \tau$ applied to them, we regard Row A as the most fundamental of the six rows. The quotient relationships are fully manifested only within Column 0, and this column has the nicest set of answers amongst the columns. This column can be used to guide the application of the three symmetries to form the assemblages of the other five rows from those of Row A. Given the summation process cited above, we regard Column 0 as being the most fundamental column.

In combinatorics, quotient relationships are among the most fundamental relationships between counts. The nicest of the quotient relationships presented in Section 3 is Quotient Relationship (a). It cannot be stated unless Row A is present!

Once Row A and Column 0 have been included, it would not make sense to exclude Row D. Answer D3 is a nice formula for the interesting count of stacks of books. Counts D3 and A3 satisfy Quotient Relationship (b). Column 4 brings in the famous counting quantities $B(m)$, $p(m)$, and $2^{m-1}$. The one parameter counts of Column 4 are more fundamental than the two parameter counts of Columns 2 and 3. This column of one parameter counts led to the definitions of Columns 5-10 and the mutually beneficial interactions with the OEIS described in Section 7. "Clicking around" within the OEIS using the on-line version of Table 3 further indicates the value of including Rows A and D and Column 4, and the organizing utility of the tabular viewpoint.

**ACKNOWLEDGMENTS.** I am indebted to Doug West, Rob Donnelly, and Tom Brylawski for essential extensive conversations. I would also like to thank Paul Proctor, Bruce Sagan, Alina Badus, Scott Lewis, and Cheryl Gann for some helpful comments.

**Table 2: Thirty Assemblage Counting Answers**

| Number of ways to group m Items. | Nature of Items | Nature of Batches | (0) $a_i$ Batches Contain i Items | (1) b Blips (# of items ≤ 1) | (2) b Batches (any # of items) | (3) b Blocks (# of items ≥ 1) | (4) Any Number of Blocks |
|---|---|---|---|---|---|---|---|
| **Part I** Sequences of Batches (A) | Disting'ble | Lists | $\binom{b}{\alpha} m!$ | $b_{(m)}$ | $b^{(m)}$ | $\binom{m-1}{b-1} m!$ | $2^{m-1} m!$ |
| (B) | Disting'ble | Sets | $\binom{b}{\alpha}\binom{m}{i}$ | $b_{(m)}$ | $b^m$ | $b!\,S(m,b)$ | $\sum_{b=1}^{m} b!\,S(m,b)$ |
| (C) | Identical | Bunches | $\binom{b}{\alpha}$ | $\binom{b}{m} := \dfrac{b_{(m)}}{m!}$ | $\left(\!\binom{b}{m}\!\right) := \dfrac{b^{(m)}}{m!}$ | $\binom{m-1}{b-1}$ | $2^{m-1}$ |
| **Part II** Collections of Batches (D) | Disting'ble | Lists | $\dfrac{m!}{\alpha^{+}!}$ | $\begin{cases} 1 & \text{if } b \ge m \\ 0 & \text{if } b < m \end{cases}$ | $\sum_{d=1}^{b} \dfrac{m!}{d!}\binom{m-1}{d-1}$ | $\dfrac{m!}{b!}\binom{m-1}{b-1}$ | $\sum_{b=1}^{m} \dfrac{m!}{b!}\binom{m-1}{b-1}$ |
| (E) | Disting'ble | Sets | $\dfrac{1}{\alpha^{+}!}\binom{m}{i}$ | $\begin{cases} 1 & \text{if } b \ge m \\ 0 & \text{if } b < m \end{cases}$ | $\sum_{d=1}^{b} S(m,d)$ | $S(m,b)$ | $B(m)$ |
| (F) | Identical | Bunches | $1$ | $\begin{cases} 1 & \text{if } b \ge m \\ 0 & \text{if } b < m \end{cases}$ | $\sum_{d=1}^{b} p(m,d)$ | $p(m,b)$ | $p(m)$ |